\documentclass{article}
\usepackage{amsmath}
\usepackage{amsfonts}

\input{tcilatex}
\begin{document}

\begin{center}
\textbf{ROTATIONAL SURFACES IN ISOTROPIC SPACES SATISFYING WEINGARTEN
CONDITIONS}

\bigskip

\textbf{Alper Osman \"{O}\u{g}renmi\c{s}}

Department of Mathematics, Faculty of Science

Firat University, Elazig, 23119, Turkey

aogrenmis@firat.edu.tr

\bigskip
\end{center}

\textbf{Abstract. }In this paper, we study the rotational surfaces in the
isotropic 3-space $\mathbb{I}^{3}$ satisfying Weingarten conditions in terms
of the relative curvature $K$ (analogue of the Gaussian curvature) and the
isotropic mean curvature $H$. In particular, we classify such surfaces of
linear Weingarten type in $\mathbb{I}^{3}.$

\bigskip

\textbf{Keywords:} Isotropic space; rotational surface; Weingarten surface.

\textbf{Math. Subject Classification 2010: }$53A35$, $53A40$, $53B25$.

\section{Introduction}

The work of surfaces with special properties in the isotropic 3-space $%
\mathbb{I}^{3}$ has important applications in several applied sciences,
e.g., computer science, Image Processing, architectural design and
microeconomics, see \cite{3,4,6,8}, \cite{29}-\cite{31}.

Differential geometry of isotropic spaces have been introduced by K.
Strubecker \cite{37}, H. Sachs \cite{32}-\cite{34}, D. Palman \cite{27} and
others.

I. Kamenarovic (\cite{17,18}), B. Pavkovic (\cite{28}), Z. M. Sipus (\cite%
{35}) and M.E. Aydin (\cite{1,2}) have studied some classes of surfaces in $%
\mathbb{I}^{3}.$

On the other hand, let $\mathcal{M}$ be a regular surface of a Euclidean
3-space $\mathbb{%
\mathbb{R}
}^{3}$. For general references on the geometry of surfaces see \cite{12,15}.

Denote $\nabla $ the Levi-Civita connection of $\mathbb{%
\mathbb{R}
}^{3}$ and $N$ the normal vector field to $\mathcal{M}.$ Then the operator
given by%
\begin{equation*}
S\left( v\right) =-\nabla _{v}N,
\end{equation*}%
is called the \textit{shape operator}, where $v$ is a tangent vector field
to $\mathcal{M}.$ It measures how $\mathcal{M}$ bends in different
directions. The eigenvalues of $S$ are called the \textit{principal
curvatures} donoted by $\kappa _{1}$ and $\kappa _{2}.$

The arithmetic mean of the principal curvatures are called the \textit{mean
curvature}, $H=\frac{1}{2}\left( \kappa _{1}+\kappa _{2}\right) .$ The 
\textit{Gaussian curvature} is defined by $K=\kappa _{1}\kappa _{2}.$

A surface $\mathcal{M}$ in $%
\mathbb{R}
^{3}$ is called a \textit{Weingarten surface (W-surface)} if it satisfies
the following non-trivial functional relation 
\begin{equation*}
\phi \left( \kappa _{1},\kappa _{2}\right) =0
\end{equation*}%
for a smooth function $\phi $ of two variables. The above relation implies
the following%
\begin{equation*}
\delta \left( K,H\right) =0,
\end{equation*}%
which is the equivalent to the vanishing of the corresponding Jacobian
determinant, i.e. $\left\vert \frac{\partial \left( K,H\right) }{\partial
\left( u,v\right) }\right\vert =0$ for a coordinate pair $\left( u,v\right) $
on $\mathcal{M}.$

If $\mathcal{M}$ fulfills the following condition 
\begin{equation*}
c_{1}H+c_{2}K=c_{3},\text{ }c_{i}\in 
\mathbb{R}
,\text{ }\left( c_{1},c_{2},c_{3}\right) \neq \left( 0,0,0\right) ,\text{ }%
i=1,2,3,
\end{equation*}%
then it is called a\textit{\ linear Weingarten surface (LW-surface). }In the
particular case $c_{1}=0$ (resp. $c_{2}=0$)$,$ the LW-surfaces are indeed
the surfaces with constant Gaussian curvature (resp. mean curvature). These
phenomenal surfaces have been stuied by many geometers in various ambient
spaces, see \cite{14,20}, \cite{22}-\cite{24}, \cite{26}, \cite{38}.

The motivation of the present paper is to study Weingarten surfaces, in
particular Weingarten rotational surfaces, in the isotropic 3-space $\mathbb{%
I}^{3}$ which is one the Cayley--Klein spaces.

Most recently, M.E. Aydin (\cite{2}) classified the helicoidal surfaces in $%
\mathbb{I}^{3}$, which are natural generalization of the rotational
surfaces, with constant curvature and analyzed some special curves on such
surfaces.

In the present paper, we provide that the rotational surfaces in $\mathbb{I}%
^{3}$ are evidently Weingarten ones. Then we classified LW-rotational
surfaces in $\mathbb{I}^{3}$ satisfying the following relation%
\begin{equation*}
K=m_{0}H+n_{0},\text{ }m_{0},n_{0}\in 
\mathbb{R}
,
\end{equation*}%
in which $K$ is the relative curvature and $H$ isotropic mean curvature.

\section{Preliminaries}

The isotropic 3-space $\mathbb{I}^{3}$ is obtained from the 3-dimensional
projective space $P\left( 
\mathbb{R}
^{3}\right) $ with the absolute figure which is an ordered triple $\left(
p,l_{1},l_{2}\right) $, where $p$ is a plane in $P\left( 
\mathbb{R}
^{3}\right) $ and $l_{1},l_{2}$ are two complex-conjugate straight lines in $%
p$ (see \cite{35}). The homogeneous coordinates in $P\left( 
\mathbb{R}
^{3}\right) $ are introduced in such a way that the \textit{absolute plane} $%
p$ is given by $x_{0}=0$ and the \textit{absolute lines} $l_{1},l_{2}$ by $%
x_{0}=x_{1}+ix_{2}=0,$ $x_{0}=x_{1}-ix_{2}=0.$ The intersection point $%
P(0:0:0:1)$ of these two lines is called the \textit{absolute point}. The
group of motions of $\mathbb{I}^{3}$ is a six-parameter group given in the
normal form (in affine coordinates) $\mathbf{x}=\frac{x_{1}}{x_{0}},$ $%
\mathbf{y}=\frac{x_{2}}{x_{0}},$\textbf{\ }$\mathbf{z}=\frac{x_{3}}{x_{0}}$
by

\begin{equation}
\left( \mathbf{x},\mathbf{y},\mathbf{z}\right) \longmapsto \left( \mathbf{x}%
^{\prime },\mathbf{y}^{\prime },\mathbf{z}^{\prime }\right) :\left\{ 
\begin{array}{l}
\mathbf{x}^{\prime }=c_{1}+\mathbf{x}\cos c_{2}-\mathbf{y}\sin c_{2}, \\ 
\mathbf{y}^{\prime }=c_{3}+\mathbf{x}\sin c_{2}+\mathbf{y}\cos c_{2}, \\ 
\mathbf{z}^{\prime }=c_{4}+c_{5}\mathbf{x}+c_{6}\mathbf{y}+\mathbf{z},%
\end{array}%
\right.  \tag{2.1}
\end{equation}%
for $c_{1},...,c_{6}\in 
\mathbb{R}
.$

Such affine transformations are called \textit{isotropic congruence
transformations }or \textit{i-motions. }

Consider the points $p_{1}=\left( x_{1},x_{2},x_{3}\right) $ and $%
p_{2}=\left( y_{1},y_{2},y_{3}\right) .$ The \textit{isotropic distance},
so-called \textit{i-distance }of two points $p_{1}$ and $p_{2}$ is defined by%
\begin{equation*}
\left\Vert p_{1}-p_{2}\right\Vert _{i}=\left( \left( y_{1}-x_{1}\right)
^{2}+\left( y_{2}-x_{2}\right) ^{2}\right) ^{\frac{1}{2}}.
\end{equation*}%
The i-metric is degenerate along the lines in $\mathbf{z}-$direction, and
such lines are called \textit{isotropic} lines.

\textit{Planes, circles and spheres}. There are two types of planes in $%
\mathbb{I}^{3}$ (\cite{29}-\cite{31}).

(1) \textit{Non-isotropic planes} are planes non-parallel to the $\mathbf{z}%
- $direction. In these planes we basically have an Euclidean metric: This is
not the one we are used to, since we have to make the usual Euclidean
measurements in the top view. An \textit{i-circle} (of \textit{elliptic type}%
) in a non-isotropic plane $p$ is an ellipse, whose top view is an Euclidean
circle. Such an i-circle with center $c_{0}\in p$ and radius $r$ is the set
of all points\textbf{\ }$x\in p$ with $\left\Vert x-c_{0}\right\Vert _{i}=r.$

(2) \textit{Isotropic planes }are planes parallel to the $\mathbf{z}-$axis.
There, $\mathbb{I}^{3}$ induces an isotropic metric. An \textit{i-circle }%
(of \textit{parabolic type}) is a parabola with $\mathbf{z}-$parallel axis
and thus it lies in an isotropic plane

An i-circle of parabolic type is not the iso-distance set of a fixed point,
but it may be seen as a curve with constant isotropic curvature: A curve $%
\alpha $ in an isotropic plane P (without loss of generality we set $P:%
\mathbf{y}=0$) which does not possess isotropic tangents can be written as
graph $\mathbf{z}=f(\mathbf{x})$. Then, the \textit{i-curvature} of $\alpha $
at $\mathbf{x}=s_{0}$ is given by the second derivative $\kappa _{i}\left(
s_{0}\right) =f^{\prime \prime }\left( s_{0}\right) $. For an i-circle of
parabolic type $f$ is quadratic and thus $\kappa _{i}$ is constant.

There are also two types of \textit{isotropic spheres}. An \textit{i-sphere
of the cylindrical type} is the set of all points $x\in \mathbb{I}^{3}$ with 
$\left\Vert x-c_{0}\right\Vert _{i}=r$. Speaking in an Euclidean way, such a
sphere is a right circular cylinder with $\mathbf{z}\mathit{-}$parallel
rulings; its top view is the Euclidean circle with center $c_{0}$ and radius 
$r$. The more interesting and important type of spheres are the \textit{%
i-spheres of parabolic type,}%
\begin{equation*}
\mathbf{z}=\frac{A}{2}\left( \mathbf{x}^{2}+\mathbf{y}^{2}\right) +B\mathbf{x%
}+C\mathbf{y}+D,\text{ \ \ }A\neq 0.
\end{equation*}

From an Euclidean perspective, they are paraboloids of revolution with $%
\mathbf{z-}$parallel axis. The intersections of these i-spheres with planes $%
p$ are i-circles. If $p$ is non-isotropic, then the intersection is an
i-circle of elliptic type. If $p$ is isotropic, the intersection curve is an
i-circle of parabolic type.

\textit{Curvature theory of surfaces. }A surface $\mathcal{M}$ immersed in $%
\mathbb{I}^{3}$ is called \textit{admissible} if it has no isotropic tangent
planes. We restrict our framework to admissible regular surfaces. For such a
surface $\mathcal{M}$, the coefficients $E,F,G$ of its first fundamental
form are calculated with respect to the induced metric.

The normal field of $\mathcal{M}$ is always the isotropic vector$.$ The
coefficients $L,M,N$ of the second fundamental form of $\mathcal{M}$ are
calculated with respect to the normal field of $\mathcal{M}$ (for details,
see \cite{33}, p. 155).

The \textit{relative curvature }(so called \textit{isotropic Gaussian
curvature}) and \textit{isotropic mean curvature} are defined by%
\begin{equation}
K=\frac{LM-N^{2}}{EG-F^{2}},\text{ \ \ }H=\frac{EN-2FM+GL}{2EG-F^{2}}. 
\tag{2.2}
\end{equation}%
\newline

\section{LW-rotational surfaces in $\mathbb{I}^{3}$}

Let us consider the i-motions given by $\left( 2.1\right) ,$ then the
Euclidean rotations in the isotropic space $\mathbb{I}^{3}$ is given by in
affine coordinates%
\begin{equation*}
\left\{ 
\begin{array}{l}
\mathbf{x}^{\prime }=c_{1}+\mathbf{x}\cos c_{2}-\mathbf{y}\sin c_{2}, \\ 
\mathbf{y}^{\prime }=c_{3}+\mathbf{x}\sin c_{2}+\mathbf{y}\cos c_{2}, \\ 
\mathbf{z}^{\prime }=\mathbf{z},%
\end{array}%
\right.
\end{equation*}%
where $c_{i}\in 
\mathbb{R}
.$

\bigskip

\textbf{Definition 3.1. }\textit{Let }$\alpha $\textit{\ be a curve lying in
the isotropic }$xz-$\textit{plane given by }$c\left( u\right) =\left(
u,0,g\left( u\right) \right) $\textit{\ where }$g\in C^{2},$ $\frac{dg}{du}%
\neq 0$\textit{. By rotating the curve }$c$\textit{\ around }$z-$\textit{%
axis, we obtain that the rotational surface in }$\mathbb{I}^{3}$ \textit{is
of the form}

\begin{equation}
\mathbf{X}\left( u,v\right) =\left( u\cos v,u\sin v,g\left( u\right) \right)
.  \tag{3.1}
\end{equation}%
\textit{Similarly when the profile curve }$\alpha $\textit{\ lies in the
isotropic }$yz-$\textit{plane, then the\ rotational surface in} $\mathbb{I}%
^{3}$\textit{\ is given by}

\begin{equation}
\mathbf{X}\left( u,v\right) =\left( -u\sin v,u\cos v,g\left( u\right)
\right) .  \tag{3.2}
\end{equation}

\bigskip

\textbf{Remark 3.1. }The rotational surfaces given by $\left( 3.1\right) $
and $\left( 3.2\right) $ are locally isometric and thus we only consider the
ones of the form $\left( 3.1\right) $.

\bigskip

Let $\mathcal{M}$ be the rotational surface given by $\left( 3.1\right) $ in%
\textit{\ }$\mathbb{I}^{3}.$ Then the nonzero components of first
fundamental form of $\mathcal{M}$ are calculated by induced metric from $%
\mathbb{I}^{3}$ as follows%
\begin{equation}
E=1\text{, }G=u^{2}.  \tag{3.3}
\end{equation}%
The nonzero components of second fundamental form of $\mathcal{M}$ are%
\begin{equation}
L=g^{\prime \prime },\text{ }N=ug^{\prime },  \tag{3.4}
\end{equation}%
where $g^{\prime }=\frac{dg}{du}$ and $g^{\prime \prime }=\frac{d^{2}g}{%
du^{2}}.$ From $\left( 2.2\right) ,$ $\left( 3.3\right) $ and $\left(
3.4\right) ,$ we get%
\begin{equation}
K=\frac{1}{u}g^{\prime }g^{\prime \prime },\text{ }H=\frac{1}{u}g^{\prime
}+g^{\prime \prime },  \tag{3.5}
\end{equation}%
which yields that the curvatures $K$ and $H$ depend only on the variable $u$%
, namely $\left\vert \frac{\partial \left( K,H\right) }{\partial \left(
u,v\right) }\right\vert =0.$ In the sequel, we have the following result.

\bigskip

\textbf{Theorem 3.1. }\textit{Rotational surfaces in }$\mathbb{I}^{3}$ 
\textit{are Weingarten surfaces.}

\bigskip

We are also able to investigate the LW-rotational surfaces in\textit{\ }$%
\mathbb{I}^{3}$ with the relation%
\begin{equation}
K=m_{0}H+n_{0},\text{ }m_{0},n_{0}\in 
\mathbb{R}
.  \tag{3.6}
\end{equation}%
If $m_{0}=0$ in $\left( 3.6\right) ,$ then those reduce to ones with
constant relative curvature. Thus we aim to obtain the LW-rotational
surfaces in\textit{\ }$\mathbb{I}^{3}$ with $m_{0}\neq 0.$

The following result classifies the LW-rotational surfaces satisfying $%
\left( 3.6\right) .$

\bigskip

\textbf{Theorem 3.2. }\textit{Let }$\mathcal{M}$ \textit{be a LW-rotational
surface in }$\mathbb{I}^{3}$. \textit{Then one of the following holds}

\textit{(i) }$\mathcal{M}$\textit{\ is of the form }%
\begin{equation*}
\left\{ 
\begin{array}{l}
\mathbf{X}\left( u,v\right) =\left( u\cos v,u\sin v,g\left( u\right) \right)
, \\ 
g\left( u\right) =\frac{m_{0}}{2}u^{2}\pm \frac{u}{2}\sqrt{%
c_{1}+m_{0}^{2}u^{2}}\pm c_{2}\ln \left\vert 2m_{0}\left( m_{0}u+\sqrt{%
c_{1}+m_{0}^{2}u^{2}}\right) \right\vert , \\ 
c_{1},c_{2}\in 
\mathbb{R}
\backslash \left\{ 0\right\} ;%
\end{array}%
\right.
\end{equation*}

\textit{(ii) }$\mathcal{M}$ \textit{is an elliptic paraboloid from the
Euclidean perspective, i.e.}%
\begin{equation*}
\left\{ 
\begin{array}{l}
\mathbf{X}\left( u,v\right) =\left( u\cos v,u\sin v,g\left( u\right) \right)
, \\ 
g\left( u\right) =\frac{m_{0}}{2}u^{2}+c_{3},\text{ }c_{3}\in 
\mathbb{R}
;%
\end{array}%
\right.
\end{equation*}

\textit{(iii) }$\mathcal{M}$\textit{\ is given by }%
\begin{equation*}
\left\{ 
\begin{array}{l}
\mathbf{X}\left( u,v\right) =\left( u\cos v,u\sin v,g\left( u\right) \right)
, \\ 
g\left( u\right) =\frac{m_{0}}{2}u^{2}\pm \frac{u}{2}\sqrt{c_{1}+\left(
m_{0}^{2}+n_{0}\right) u^{2}}\pm \\ 
\pm \frac{c_{1}}{m_{0}^{2}+n_{0}}\ln \left\vert 2\left( \left(
m_{0}^{2}+n_{0}\right) u+\sqrt{m_{0}^{2}+n_{0}}\sqrt{c_{1}+\left(
m_{0}^{2}+n_{0}\right) u^{2}}\right) \right\vert , \\ 
c_{1}\in 
\mathbb{R}
,\text{ }c_{1}<0.%
\end{array}%
\right.
\end{equation*}

\textbf{Proof. }Assume $\mathcal{M}$ is a LW-rotational surface in\textit{\ }%
$\mathbb{I}^{3}$ having the relation $\left( 3.6\right) .$ Then, from $%
\left( 3.5\right) ,$ it follows%
\begin{equation}
\frac{1}{u}g^{\prime }g^{\prime \prime }=m_{0}\frac{g^{\prime }+g^{\prime
\prime }}{u}+n_{0}.  \tag{3.7}
\end{equation}%
We have two cases:

\bigskip

\textbf{Case a. }$n_{0}=0.$ Hence we can rewrite $\left( 3.7\right) $ as%
\begin{equation}
g^{\prime \prime }\left( g^{\prime }-m_{0}u\right) -m_{0}g^{\prime }=0. 
\tag{3.8}
\end{equation}%
If $g^{\prime }=m_{0}u$ in $\left( 3.8\right) ,$ then $g^{\prime }$ and $%
m_{0}$ vanish which is not possible. Then we have%
\begin{equation}
g^{\prime \prime }-\frac{m_{0}g^{\prime }}{g^{\prime }-m_{0}u}=0.  \tag{3.9}
\end{equation}%
By solving $\left( 3.9\right) ,$ we obtain 
\begin{equation*}
g\left( u\right) =\frac{m_{0}}{2}u^{2}\pm \frac{u}{2}\sqrt{%
e^{2c_{1}}+m_{0}^{2}u^{2}}\pm \frac{e^{2c_{1}}}{2m_{0}}\ln \left\vert
2m_{0}\left( m_{0}u+\sqrt{e^{2c_{1}}+m_{0}^{2}u^{2}}\right) \right\vert ,
\end{equation*}%
$c_{1}\in 
\mathbb{R}
,$ which gives the statement (i) of the theorem.

\bigskip

\textbf{Case b.} $n_{0}\neq 0.$ Then we have from $\left( 3.7\right) $%
\begin{equation}
g^{\prime \prime }\left( g^{\prime }-m_{0}u\right) -m_{0}g^{\prime }=n_{0}u.
\tag{3.10}
\end{equation}%
When $g^{\prime }=m_{0}u,$ then $g\left( u\right) =\frac{m_{0}}{2}%
u^{2}+c_{2},$ $c_{2}\in 
\mathbb{R}
$ and $n_{0}=-m_{0}^{2}.$ This implies the statement (ii) of the theorem.

Otherwise, we conclude from $\left( 3.10\right) $ that%
\begin{equation}
g^{\prime \prime }-\frac{m_{0}g^{\prime }}{g^{\prime }-m_{0}u}=\frac{n_{0}u}{%
g^{\prime }-m_{0}u}.  \tag{3.11}
\end{equation}%
After solving $\left( 3.11\right) ,$ we derive%
\begin{equation*}
\begin{array}{l}
g\left( u\right) =\frac{m_{0}}{2}u^{2}\pm \frac{u}{2}\sqrt{%
-e^{2c_{3}}+\left( m_{0}^{2}+n_{0}\right) u^{2}}\mp \\ 
\pm \frac{e^{2c_{3}}}{m_{0}^{2}+n_{0}}\ln \left\vert 2\left( \left(
m_{0}^{2}+n_{0}\right) u+\sqrt{m_{0}^{2}+n_{0}}\sqrt{-e^{2c_{3}}+\left(
m_{0}^{2}+n_{0}\right) u^{2}}\right) \right\vert ,%
\end{array}%
\end{equation*}%
$c_{3}\in 
\mathbb{R}
.$ Therefore the proof is completed.

\bigskip

\textbf{Example 3.1. }\textit{Consider the elliptic paraboloid in }$\mathbb{I%
}^{3}$\textit{\ from the Euclidean perspective given by}%
\begin{equation*}
\mathbf{X}\left( u,v\right) =\left( u\cos v,u\sin v,0.25u^{2}\right) ,\text{ 
}\left( u,v\right) \in \left[ 0,2\pi \right] .
\end{equation*}%
Then $K=0.25,$ $H=1,$ $m_{0}=0.5$ and $n_{0}=-0.25.$ We plot it as in Fig. 1.%
\begin{gather*}
\FRAME{itbpF}{1.6734in}{1.4909in}{0in}{}{}{Figure}{\special{language
"Scientific Word";type "GRAPHIC";maintain-aspect-ratio TRUE;display
"USEDEF";valid_file "T";width 1.6734in;height 1.4909in;depth
0in;original-width 1.6362in;original-height 1.4555in;cropleft "0";croptop
"1";cropright "1";cropbottom "0";tempfilename
'O4Y09I00.wmf';tempfile-properties "XPR";}} \\
\mathbf{Fig}\text{ }\mathbf{1.}\text{ \textit{LW-rotational surface with }}%
m_{0}=0.5,\text{ }n_{0}=-0.25
\end{gather*}

\section{Rotational surfaces in $\mathbb{I}^{3}$ with $H/K=const.$}

The authors in \cite{7} introduced a new kind of curvature for the
hypersurfaces of Euclidean $n-$spaces, called by amalgamatic curvature and
explored its geometric meaning by proving an inequality related to the
absolute mean curvature of the hypersurface. In the particular case $n=3$,
the amalgamatic curvature is indeed the harmonic ratio of the principal
curvatures of any given surface, i.e., the ratio of the Gaussian curvature
and the mean curvature.

By considering this argument, we can consider the rotational surfaces in $%
\mathbb{I}^{3}$ satisfying $H/K=const.$ Thus the statement (i) of Theorem
3.2 is indeed a classification of the rotational surfaces in $\mathbb{I}^{3}$
satisfying $H/K=const.$

Therefore, we have the following trivial result.

\bigskip

\textbf{Corollary 4.1. }\textit{Let }$\mathcal{M}$ \textit{be a
LW-rotational surface in }$\mathbb{I}^{3}$ \textit{satisfying} $H/K=\frac{1}{%
m_{0}},$ $m_{0}\in 
\mathbb{R}
\backslash \left\{ 0\right\} .$ \textit{Then it is of the form}%
\begin{equation}
\left\{ 
\begin{array}{l}
\mathbf{X}\left( u,v\right) =\left( u\cos v,u\sin v,g\left( u\right) \right)
, \\ 
g\left( u\right) =\frac{m_{0}}{2}u^{2}\pm \frac{u}{2}\sqrt{%
c_{1}+m_{0}^{2}u^{2}}\pm c_{2}\ln \left\vert 2m_{0}\left( m_{0}u+\sqrt{%
c_{1}+m_{0}^{2}u^{2}}\right) \right\vert , \\ 
c_{1},c_{2}\in 
\mathbb{R}
\backslash \left\{ 0\right\} ;%
\end{array}%
\right.  \tag{4.1}
\end{equation}

\bigskip

\textbf{Example 4.1. }\textit{Take }$\lambda _{0}=0.5$\textit{\ and }$%
c_{1}=\ln 2$\textit{\ in }$\left( 4.1\right) .$\textit{\ Then we obtain a
rotational surface in }$\mathbb{I}^{3}$ \textit{with }$H/K=1$ \textit{given
by}%
\begin{equation*}
\mathbf{X}\left( u,v\right) =\left( u\cos v,u\sin v,u^{2}+u\sqrt{1+u^{2}}%
+\ln \left\vert 2\left( u+\sqrt{1+u^{2}}\right) \right\vert \right) ,
\end{equation*}%
where $u\in \left[ 0,2\pi \right] $, $v\in \left[ 0,\pi /2\right] .$ Then it
can be plotted as in Fig. 2.%
\begin{gather*}
\FRAME{itbpF}{1.0309in}{1.8766in}{0in}{}{}{Figure}{\special{language
"Scientific Word";type "GRAPHIC";maintain-aspect-ratio TRUE;display
"USEDEF";valid_file "T";width 1.0309in;height 1.8766in;depth
0in;original-width 0.998in;original-height 1.8395in;cropleft "0";croptop
"1";cropright "1";cropbottom "0";tempfilename
'O4Y09W01.wmf';tempfile-properties "XPR";}} \\
\mathbf{Fig}\text{ }\mathbf{2.}\text{ \textit{Rotational surface with }}H/K=1
\end{gather*}


\begin{thebibliography}{99}
\bibitem{1} M.E. Aydin, \textit{A generalization of translation surfaces
with constant curvature in the isotropic space, }J. Geom., 2015, DOI
10.1007/s00022-015-0292-0.

\bibitem{2} M.E. Aydin, \textit{Classification results on surfaces in the
isotropic 3-space}, arXiv:1601.03190v1 [math.DG], 2016.

\bibitem{3} M.E. Aydin and A. Mihai, \textit{Classification of quasi-sum
production functions with Allen determinants}, Filomat \textbf{29(6)}
(2015), 1351--1359.

\bibitem{4} M.E. Aydin and A. Mihai, \textit{Translation hypersurfaces and
Tzitzeica translation hypersurfaces of the Euclidean space}, Proc. Ro. Acad.
Series A \textbf{16(4) }(2015), 477-483.

\bibitem{5} C. Baikoussis and T. Koufogioros, \textit{Helicoidal surface
with prescribed mean or Gauss curvature}, J. Geom. \textbf{63} (1998),
25--29.

\bibitem{6} B. Y. Chen, S. Decu and L. Verstraelen, \textit{Notes on
isotropic geometry of production models}, Kragujevac J. Math. \textbf{37(2)}
(2013), 217--220.

\bibitem{7} C. T. R. Conley, R. Etnyre, B. Gardener, L. H. Odom and B. D.
Suceava, \textit{New curvature inequalities for hypersurfaces in the
Euclidean ambient space}, Taiwanese J. Math. 17(3) (2013), 885--895.

\bibitem{8} S. Decu, L. Verstraelen, \textit{A note on the isotropical
geometry of production surfaces}, Kragujevac J. Math. \textbf{38(1)\ }%
(2014), 23--33.

\bibitem{9} G. Delaunay,\textit{\ Sur la surface de revolution dont la
courbure moyenne est constante}, J. Math. Pures Appl. Series\textbf{\ 6(1)}
(1841), 309-320.

\bibitem{10} F. Dillen and W. Kuhnel, \textit{Ruled Weingarten surfaces in
Minkowski 3-space}, Manuscripta Math., \textbf{98} (1999), 307-320.

\bibitem{11} M.P. Do Carmo and M. Dajczer,\textit{\ Helicoidal surfaces with
constant mean curvature}, Tohoku Math. J. \textbf{34 }(1982), 425-435.

\bibitem{12} M.P. Do Carmo, Differential geometry of curves and surfaces,
Prentice Hall: Englewood Cliffs, NJ, 1976.

\bibitem{13} Z. Erjavec, B. Divjak and D. Horvat, \textit{The general
solutions of Frenet's system in the equiform geometry of the Galilean,
pseudo-Galilean, simple isotropic and double isotropic space, }Int. Math.
Forum \textbf{6(17) }(2011), 837-856.

\bibitem{14} J. A. Galvez, A. Martinez and F. Milan, \textit{Linear
Weingarten surfaces in} $R^{3}$, Monatsh. Math., \textbf{138 }(2003),
133-144.

\bibitem{15} A. Gray, Modern differential geometry of curves and surfaces
with mathematica. CRC Press LLC, 1998.

\bibitem{16} Z.H. Hou and F. Ji, \textit{Helicoidal surfaces with} $H^{2}=K$ 
\textit{in Minkowski 3-space}, J. Math. Anal. Appl.\textbf{\ 325} (2007),
101--113.

\bibitem{17} I. Kamenarovic,\textit{\ On line complexes in the isotropic
space} $I_{3}^{(1)},$ Glasnik Matematicki \textbf{17(37)} (1982), 321-329.

\bibitem{18} I. Kamenarovic, \textit{Associated curves on ruled surfaces in
the isotropic space} $I_{3}^{(1)},$ Glasnik Matematicki \textbf{29(49)}
(1994), 363-370.

\bibitem{19} K. Kenmotsu, \textit{Surfaces of revolution with prescribed
mean curvature}, Tohoku Math. J. 32 (1980), 147-153.

\bibitem{20} M. H Kim and D. W. Yoon,\textit{\ Weingarten quadric surfaces
in a Euclidean 3-space}, Turk. J. Math. \textbf{35} (2011), 479-485.

\bibitem{21} J. J. Koenderink and A. van Doorn,\textit{\ Image processing
done right}, Lecture Notes in Computer Science \textbf{2350} (2002),
158--172.

\bibitem{22} W. Kuhnel, \textit{Ruled W-surfaces}, Arch. Math. \textbf{62 }%
(1994), 475-480.

\bibitem{23} C.W. Lee, \textit{Linear Weingarten rotational surfaces in
pseudo-Galilean 3-space}, Int. J. Math. Anal.\textbf{\ 9(50) }(2015), 2469 -
2483.

\bibitem{24} H. Liu and G. Liu, \textit{Weingarten rotation surfaces in
3-dimensional de Sitter space, }J. Geom. \textbf{79} (2004), 156 -- 168.

\bibitem{25} R. Lopez and E. Demir,\textit{\ Helicoidal surfaces in
Minkowski space with constant mean curvature and constant Gauss curvature},
Cent. Eur. J. Math. \textbf{12(9)} (2014), 1349-1361.

\bibitem{26} R. Lopez, \textit{Rotational linear Weingarten surfaces of
hyperbolic type}, Israel J. Math. \textbf{167 }(2008), 283--301.

\bibitem{27} D. Palman, \textit{Spharische quartiken auf dem torus im
einfach isotropen raum}, Glasnik Matematicki \textbf{14(34)} (1979), 345-357.

\bibitem{28} B. Pavkovic, \textit{An interpretation of the relative
curvatures for surfaces in the isotropic space}, Glasnik Matematicki \textbf{%
15(35)} (1980), 149-152.

\bibitem{29} H. Pottmann and K. Opitz, \textit{Curvature analysis and
visualization for functions defined on Euclidean spaces or surfaces},
Comput. Aided Geom. Design \textbf{11} (1994), 655--674.

\bibitem{30} H. Pottmann and Y. Liu, \textit{Discrete surfaces of isotropic
geometry with applications in architecture}. In: Martin, R., Sabin, M.,
Winkler, J. (eds.) The Mathematics of Surfaces, pp. 341--363. Lecture Notes
in Computer Science 4647. Springer (2007).

\bibitem{31} H. Pottmann, P. Grohs and N.J. Mitra, \textit{Laguerre minimal
surfaces, isotropic geometry and linear elasticity}, Adv. Comput. Math. 
\textbf{31} (2009), 391--419

\bibitem{32} H. Sachs, Ebene Isotrope Geometrie, Vieweg-Verlag,
Braunschweig, Wiesbaden, 1990.

\bibitem{33} H. Sachs, Isotrope Geometrie des Raumes, Vieweg Verlag,
Braunschweig, 1990.

\bibitem{34} H. Sachs, \textit{Zur Geometrie der Hyperspharen in
n-dimensionalen einfach isotropen Raum,} Jour. f. d. reine u. angew. Math. 
\textbf{298} (1978), 199-217.

\bibitem{35} Z. M. Sipus, \textit{Translation surfaces of constant
curvatures in a simply isotropic space}, Period. Math. Hung. \textbf{68 }%
(2014), 160--175

\bibitem{36} Z. M. Sipus and B. Divjak, \textit{Curves in n-dimensional
k-isotropic space, }Glasnik Matematicki \textbf{33(53)} (1998), 267-286.

\bibitem{37} K. Strubecker, \textit{Differentialgeometrie des isotropen
Raumes III}, Flachentheorie, Math. Zeitsch. \textbf{48 }(1942), 369-427.

\bibitem{38} D. W. Yoon, Y. Tuncer and M. K. Karacan, \textit{Non-degenerate
quadric surfaces of Weingarten type}, Annales Polonici Math. \textbf{107 }%
(2013), 59-69.
\end{thebibliography}
\end{document}